\documentstyle[12pt]{article}
\newtheorem{defn}{{\sc Definition}}[section]
\newtheorem{thm}{Theorem}[section]
\newtheorem{Theorem}{Theorem}[section]
\newtheorem{lem}[thm]{Lemma}

\newtheorem{cor}[thm]{Corollary}
\newcommand{\sthm}{\begin{Theorem}}         
\newcommand{\ethm}{\end{Theorem}}           
\newtheorem{Corollary}[Theorem]{Corollary}   
\newcommand{\scor}{\begin{Corollary}}       
\newcommand{\ecor}{\end{Corollary}}         
\newcommand{\pf}{ \par \vspace{1ex} \noindent {\sc Proof.} \hspace{2mm}}

\begin{document}
\title{Uniqueness of coexistence state with small perturbation}
\author{Joon Hyuk Kang}
\date{}
\maketitle
\begin{abstract}
In ~\cite{ko02b}, we established some sufficient condition for the uniqueness of the positive solution to the general elliptic system for several competing species of animals
$$\left\{\begin{array}{rll}
\Delta u^{i} + u^{i}(h_{i}(u^{i}) - g_{i}(u^{1},...,u^{i-1},u^{i+1},...,u^{N})) & = & 0\;\;\mbox{in}\;\;\Omega,\\
u^{i} & = & 0\;\;\mbox{on}\;\;\partial\Omega, 
\end{array} \right.$$
for $i = 1,...,N$.
In this paper, we try to extend the uniqueness result by perturbing the reproduction and self-limitation functions $h_{i}$'s of the above model. The techniques used in this paper are super-sub solutions, maximum principles and spectrum estimates. The arguments also rely on some detailed properties for the solution of logistic equations. 
\end{abstract}
\footnotetext{
Research supported by Andrews University Faculty Research Grant 2001}
\section{Introduction}
In the last decade, a lot of research has been focused on existence and uniqueness of steady state(positive and time independent solution) to the general competition model of several species of animals
$$\left\{\begin{array}{rll}
u^{i}_{t}(x,t) & = & \Delta u^{i}(x,t) + u^{i}(x,t)(h_{i}(u^{i}(x,t))\\
& & - g_{i}(u^{1}(x,t),...,u^{i-1}(x,t),u^{i+1}(x,t),...,u^{N}(x,t)))\;\;\mbox{in}\;\;\Omega \times R^{+},\\
u^{i} & = & 0\;\;\mbox{on}\;\;\partial\Omega,
\end{array} \right.$$
or equivalently, the positive solution to
\begin{equation}
\left\{\begin{array}{rll}
\Delta u^{i}(x) & + & u^{i}(x)(h_{i}(u^{i}(x)) - g_{i}(u^{1}(x),...,u_{i-1}(x),u^{i+1}(x),...,u^{N}(x)))\\
 & = & 0\;\;\mbox{in}\;\;\Omega,\\
u^{i} & = & 0\;\;\mbox{on}\;\;\partial\Omega. 
\end{array} \right.\label{eq:a}
\end{equation}
for $i = 1,...,N$.\\
\\[0.1in]
In ~\cite{ko02b}, we established the following uniqueness result.
\begin{thm}\label{thm:1.1} Suppose\\
$(U1)$ $g_{i},h_{i} \in C^{1}$ for $i = 1,2,...,N$,\\
$(U2)$ $h_{i},-g_{i}$ are strictly decreasing for $i = 1,2,...,N$,\\
$(U3)$ $g_{i}(0,...,0) = 0$ for $i = 1,2,...,N$,\\
$(U4)$ there are $k_{1},k_{2},...,k_{N} > 0$ such that $h_{i}(0) > \lambda_{1} + g_{i}(k_{1},...,k_{i-1},k_{i+1},...,k_{N})$ and $h_{i}(u^{i}) < 0$ for $u^{i} \geq k_{i}$, where $\lambda_{1}$ is the smallest eigenvalue of $-\Delta$ with the homogeneous boundary condition.\\
Then\\
$(A)$ (\ref{eq:a}) has a solution $(u^{1},...,u^{N})$ with 
$$\theta_{h_{i}-g_{i}(k_{1},...,k_{i-1},k_{i+1},...,k_{N})} < u^{i} < \theta_{h_{i}}$$
in $\Omega$ for $i = 1,...,N$. Conversely, any solution $(u^{1},...,u^{N})$ of (\ref{eq:a}) with $u^{i} > 0$ for all $i = 1,...,N$ in $\Omega$ must satisfy these inequalities.\\
$(B)$ If $-2\sup(h_{i}') > \sum_{j=1,j\neq i}^{N}(\sup(\frac{\partial g_{i}}{\partial x_{j}}) + K\sup(\frac{\partial g_{j}}{\partial x_{i}}))$, where \\
$K = \sup_{\Omega,i\neq j}\frac{\theta_{h_{j}}}{\theta_{h_{i}-g_{i}(k_{1},...,k_{i-1},k_{i+1},...,k_{N})}}$, then (\ref{eq:a}) has a unique coexistence state.(If we choose $\epsilon = 1$ in the proof, this is true.)  
\end{thm}
Biologically, we can interpret the condition in Theorem \ref{thm:1.1} as follows. The functions $g_{i}'s,h_{i}'s$ describe how species interact among themselves and with others. Hence, the conditions imply that species interact strongly among themselves and weakly with others. 
\\[0.1in]
The question in this paper concerns small perturbation of $h_{i}'s$ without losing the uniqueness of coexistence state of (\ref{eq:a}) when $h_{i}'s$ are nicer bounded functions. The conclusion says the $N$ species may have small relaxation with which they can still coexist peacefully. 

\section{Preliminaries}
In this section we will state some preliminary results which will be useful for our later arguments.
\begin{defn}\label{defn:2.1}(Super and sub solutions)
\begin{equation}
\left\{ \begin{array}{l}
\Delta u + f(x,u) = 0\;\;\mbox{in}\;\;\Omega,\\
u|_{\partial\Omega} = 0
\end{array} \right. \label{eq:pa}
\end{equation}
where $f \in C^{\alpha}(\bar{\Omega} \times R)$ and $\Omega$ is a bounded domain in $R^{n}$.\\
$(A)$ A function $\bar{u} \in C^{2,\alpha}(\bar{\Omega})$ satisfying 
$$\left\{ \begin{array}{l}
\Delta \bar{u} + f(x,\bar{u}) \leq 0\;\;\mbox{in}\;\;\Omega,\\
\bar{u}|_{\partial\Omega} \geq 0
\end{array} \right.$$
is called an super solution to (\ref{eq:pa}).\\
$(B)$ A function $\underline{u} \in C^{2,\alpha}(\bar{\Omega})$ satisfying 
$$\left\{ \begin{array}{l}
\Delta \underline{u} + f(x,\underline{u}) \geq 0\;\;\mbox{in}\;\;\Omega,\\
\underline{u}|_{\partial\Omega} \leq 0
\end{array} \right.$$
is called a sub solution to (\ref{eq:pa}).
\end{defn}
\begin{lem}\label{lem:2.2}
Let $f(x,\xi) \in C^{\alpha}(\bar{\Omega} \times R)$ and let $\bar{u}, 
\underline{u} \in C^{2,\alpha}(\bar{\Omega})$ be respectively, super and sub 
solutions to (\ref{eq:pa}) which satisfy $\underline{u}(x) \leq \bar{u}(x), x 
\in \bar{\Omega}$. Then (\ref{eq:pa}) has a solution $u \in 
C^{2,\alpha}(\bar{\Omega})$ with $\underline{u}(x) \leq u(x) \leq \bar{u}(x), x 
\in \bar{\Omega}$.
\end{lem}
\vskip1pc
We also need some information on the solutions of the following logistic equations.
\begin{lem}\label{lem:2.5}(in \cite{ll91})
$$\left\{ \begin{array}{l}
\Delta u + uf(u) = 0\;\; \mbox{in}\;\; \Omega,\\
u|_{\partial\Omega} = 0, u > 0,
\end{array} \right.$$ 
where $f$ is a decreasing $C^{1}$ function such that there exists $c_{0} > 0$ such that $f(u) \leq 0$ for $u \geq c_{0}$ and $\Omega$ is a bounded domain in $R^{n}$.\\
If $f(0) > \lambda_{1}$, then the above equation has a unique positive solution, where 
$\lambda_{1}$ is the first eigenvalue of $-\Delta$ with homogeneous boundary 
condition. We denote this unique positive solution as $\theta_{f}$.
\end{lem}
The main property about this positive solution is that $\theta_{f}$ is increasing as $f$ is increasing.
\\[0.2in]
Especially, for $a > \lambda_{1}$, where $\lambda_{1}$ is the first eigenvalue of $-\Delta$ with homogeneous boundary condition, the unique positive solution of
$$\left\{ \begin{array}{l}
\Delta u + u(a - u) = 0\;\; \mbox{in}\;\; \Omega,\\
u|_{\partial\Omega} = 0, u > 0,
\end{array} \right.$$
is denoted by $\omega_{a} \equiv \theta_{a-x}$.  
Hence, $\omega_{a}$ is increasing as $a > 0$ is increasing.

\section{Main Result}
We consider the model
\begin{equation}
\left\{ \begin{array}{l} 
\left.\begin{array}{l}
\Delta u^{i} + u^{i}(h_{i}(u^{i}) - g_{i}(u^{1},...,u^{i-1},u^{i+1},...,u^{N})) = 0\\ 
\end{array} \right.\;\;\mbox{in}\;\;\Omega,\\
u^{i}|_{\partial\Omega} = 0.
\end{array} \right.\label{eq:1}
\end{equation}
Here $\Omega$ is a bounded, smooth domain in $R^{n}$ and\\ 
$(P1)$ $h_{i} \in C_{B}^{1,\alpha}$, $g_{i} \in C^{1}$ for $i = 1,2,...,N$, where $C_{B}^{m,\alpha}$ in general, is the set of decreasing, bounded and continuous functions up to $m$-th order partial derivatives whose $m$-th order partial derivatives are H$\ddot{o}$lder continuous with exponent $\alpha$.\\
$(P2)$ $h_{i}(0) > \lambda_{1}$, $g_{i}$ are strictly increasing and $g_{i}(0,...,0) = 0$ for $i = 1,2,...,N$,\\
$(P3)$ there are $k_{1},k_{2},...,k_{N} > 0$ such that $h_{i}(u^{i}) < 0$ for $u^{i} \geq k_{i}$ for $i = 1,2,...,N$.
\\[0.1in]
The following is the main theorem.
\begin{thm}\label{thm:3.1} Suppose\\
$(A)$ $h_{i}(0) > \lambda_{1}(g_{i}(\theta_{h_{1}},...,\theta_{h_{i-1}},\theta_{h_{i+1}},...,\theta_{h_{N}}))$, where in general, $\lambda_{1}(q)$ is the first eigenvalue of $-\Delta + q$ with homogeneous boundary condition, denoted by simply $\lambda_{1}$ when $q \equiv 0$.\\
$(B)$ (\ref{eq:1}) has a unique coexistence state $(u^{1},...,u^{N})$,\\
$(C)$ the Frechet derivative of (\ref{eq:1}) at $(u^{1},...,u^{N})$ is invertible.\\
Then there is a neighborhood $V$ of $(h_{1},...,h_{N})$ in $(C_{B}^{1,\alpha})^{N}$ such that if\\
$(\bar{h_{1}},...,\bar{h_{N}}) \in V$, then (\ref{eq:1}) with $(h_{1},...,h_{N}) = (\bar{h_{1}},...,\bar{h_{N}})$ has a unique coexistence state.  
\end{thm}
Theorem \ref{thm:3.1} looks like the consequence of Implicit Function Theorem. But the inverse function theorem only guaranteed the uniqueness locally. Theorem \ref{thm:3.1} concluded the global uniqueness. The techniques we will use includes naturally Implicit Function Theorem and a priori estimates on solutions of (\ref{eq:1}).
\\[0.1in]
Biologically, the first condition in this theorem indicates that the rates of self-reproduction is large. The condition of invertibility of Frechet derivative also illustrates that the rates of self-limitation is relatively larger than those of competitions which will be in Theorem \ref{thm:3.3}. Then the conclusion says that small perturbation of reproduction and self-limitation rates does not affect to the existence and uniqueness of positive steady state, i.e. they can still coexist peacefully even if there is some slight change of reproduction and self-limitation rates.
\pf Since the Frechet derivative of (\ref{eq:1}) at $(u^{1},...,u^{N})$ is invertible, by the Implicit Function Theorem, there is a neighborhood $V$ of $(h_{1},...,h_{N})$ in $(C_{B}^{1,\alpha})^{N}$ and a neighborhood $W$ of $(u^{1},...,u^{N})$ in $[C_{0}^{2,\alpha}(\bar{\Omega})]^{N}$ such that for all $(\bar{h_{1}},...,\bar{h_{N}}) \in V$, there is a unique positive solution $(u_{1},...,u_{N}) \in W$ of (\ref{eq:1}). Suppose the conclusion of the theorem is false. Then there are sequences $(\alpha_{1,n},...,\alpha_{N,n},u_{1,n},...,u_{N,n}), (\alpha_{1,n},...,\alpha_{N,n},u_{1,n}^{*},...,u_{N,n}^{*})$ in $V \times [C_{0}^{2,\alpha}(\bar{\Omega})]^{N}$ such that $(u_{1,n},...,u_{N,n})$ and $(u_{1,n}^{*},...,u_{N,n}^{*})$ are the positive solutions with $(h_{1},...,h_{N}) = (\alpha_{1,n},...,\alpha_{N,n})$ and $(u_{1,n},...,u_{N,n}) \neq (u_{1,n}^{*},...,u_{N,n}^{*})$ and $(\alpha_{1,n},...,\alpha_{N,n}) \rightarrow (h_{1},...,h_{N})$. By the Schauder's boundary estimate in elliptic theory and $(A)$ of the Theorem \ref{thm:1.1}, there is a constant $c > 0$ such that
$$|u_{i,n}|_{2,\alpha} \leq c\sup_{x\in\bar{\Omega}}(u_{i,n}(x)) \leq c\sup_{x\in\bar{\Omega}}\theta_{\alpha_{i,n}}(x)$$
for all $i = 1,...,N, n = 1,2,...$.\\
But, by the convergence of $\{\alpha_{i,n}\}_{n}$ and the monotonicity of $\theta_{f}$, we conclude that $|u_{i,n}|_{2,\alpha}$ is uniformly bounded. So, there is a uniformly convergent subsequence of $\{u_{i,n}\}$, again will be denoted by $\{u_{i,n}\}$.\\
Let 
$$\left.\begin{array}{ll}
(u_{1,n},...,u_{N,n}) \rightarrow (\bar{u_{1}},...,\bar{u_{N}}) \in (C^{2,\alpha})^{N},\\
(u_{1,n}^{*},...,u_{N,n}^{*}) \rightarrow (u_{1}^{*},...,u_{N}^{*}) \in (C^{2,\alpha})^{N}.
\end{array}\right.$$  
Then $(\bar{u_{1}},...,\bar{u_{N}})$ and $(u_{1}^{*},...,u_{N}^{*})$ are solutions of (\ref{eq:1}). Claim $\bar{u_{1}} > 0,...,\bar{u_{N}} > 0, u_{1}^{*} > 0,...,u_{N}^{*} > 0.$ It is enough to show that $\bar{u_{1}},...,\bar{u_{N}}$ are not identically zero because of the Maximum Principle. Suppose not. Without loss of generality, assume $\bar{u_{1}}$ is identically zero.\\
Let $\tilde{u_{1,n}} = \frac{u_{1,n}}{\parallel u_{1,n} \parallel_{\infty}}$ for all $n\in N$. Then for $i = 2,3,...,N$
$$\left\{ \begin{array}{rll}
\Delta \tilde{u_{1,n}} + \tilde{u_{1,n}}(\alpha_{1,n}(u_{1,n}) & - & g_{1}(u_{2,n},u_{3,n},...,u_{N,n})) = 0,\\ 
\Delta u_{i,n} + u_{i,n}(\alpha_{i,n}(u_{i,n}) & - & g_{i}(u_{1,n},u_{2,n},...,u_{i-1,n},u_{i+1,n},...,u_{N,n})) = 0,
\end{array} \right. \mbox{in $\Omega.$}$$
From the elliptic theory, $\tilde{u_{1,n}} \rightarrow \tilde{u_{1}}$ and
$$\left\{ \begin{array}{rll}
\Delta \tilde{u_{1}} + \tilde{u_{1}}(h_{1}(0) & - & g_{1}(\bar{u_{2}},\bar{u_{3}},...,\bar{u_{N}})) = 0,\\
\Delta \bar{u_{i}} + \bar{u_{i}}(h_{i}(\bar{u_{i}}) & - & g_{i}(\bar{u_{1}},\bar{u_{2}},...,\bar{u_{i-1}},\bar{u_{i+1}},...,\bar{u_{N}})) = 0, i = 2,...,N,\\
\end{array} \right. \mbox{in $\Omega,$}$$
by the continuity and uniform convergence.\\
Hence, $h_{1}(0) = \lambda_{1}(g_{1}(\bar{u_{2}},\bar{u_{3}},...,\bar{u_{N}})).$\\
Let $j = 2,...,N$.\\
If $\bar{u_{j}}$ is identically zero, then $\bar{u_{j}} \equiv 0 \leq \theta_{h_{j}}$. Suppose $\bar{u_{j}}$ is not identically zero. Then since
$$\left.\begin{array}{lll}
& & \Delta\bar{u_{j}} + \bar{u_{j}}h_{j}(\bar{u_{j}})\\
& = & \Delta\bar{u_{j}} + \bar{u_{j}}(h_{j}(\bar{u_{j}}) - g_{j}(\bar{u_{1}},...,\bar{u_{j-1}},\bar{u_{j+1}},...,\bar{u_{N}})\\
& & + g_{j}(\bar{u_{1}},...,\bar{u_{j-1}},\bar{u_{j+1}},...,\bar{u_{N}}))\\
& = & \bar{u_{j}}g_{j}(\bar{u_{1}},...,\bar{u_{j-1}},\bar{u_{j+1}},...,\bar{u_{N}}) \geq 0,
\end{array}\right.$$   
$\bar{u_{j}}$ is a sub solution of
$$\left\{\begin{array}{rll}
\Delta Z + Zh_{j}(Z) & = & 0\;\;\mbox{in}\;\;\Omega,\\
Z & = & 0\;\;\mbox{on}\;\;\partial\Omega.
\end{array}\right.$$
Since any constant which is larger than $k_{j}$ is an super solution of
$$\left.\{\begin{array}{rll}
\Delta Z + Zh_{j}(Z) & = & 0\;\;\mbox{in}\;\;\Omega,\\
Z & = & 0\;\;\mbox{on}\;\;\partial\Omega,
\end{array}\right.$$
by the uniqueness of positive solution, $\bar{u_{j}} \leq \theta_{h_{j}}$.
Consequently,
$$\left.\begin{array}{lll}
h_{1}(0) & = & \lambda_{1}(g_{1}(\bar{u_{2}},\bar{u_{3}},...,\bar{u_{N}}))\\
 & \leq & \lambda_{1}(g_{1}(\theta_{h_{2}},...,\theta_{h_{N}}))
\end{array}\right.$$
by the monotonicity of $g_{1}$ and the first eigenvalue, which contradicts our assumption. Consequently, $(\bar{u_{1}},...,\bar{u_{N}})$ and $(u_{1}^{*},...,u_{N}^{*})$ are coexistence states with $(h_{1},...,h_{N})$. But, since the coexistence state in this case is unique by assumption, $(\bar{u_{1}},...,\bar{u_{N}}) = (u_{1}^{*},...,u_{N}^{*}) = (u^{1},...,u^{N})$, which contradicts the Implicit Function Theorem. 
\vskip1pc
The proof of the theorem also tells us that if one of the species becomes extinct, in other word, if one is excluded by others, then that means the reproduction rates are small, i.e. the region condition of reproduction rates $(A)$ is reasonable.
\begin{thm}\label{thm:3.2} If $(\alpha_{1,n},...,\alpha_{N,n},u_{1,n},...,u_{N,n}) \rightarrow (h_{1},...,h_{N},u^{1},...,u^{N})$ and if $u^{j} \equiv 0$ for some $j = 1,...,N$, then
$h_{j}(0) \leq \lambda_{1}(g_{j}(\theta_{h_{1}},...,\theta_{h_{j-1}},\theta_{h_{j+1}},...,\theta_{h_{N}}))$.
\end{thm}  
The condition, invertibility of Frechet derivative, in Theorem \ref{thm:3.1} is too artificial.
Now we turn out attention to get conditions to guarantee the invertibility of the Frechet derivative.
\begin{thm}\label{thm:3.3}
Suppose $(u_{1},u_{2},...,u_{N})$ is a positive solution to (\ref{eq:1}). If
$$\left.\begin{array}{lll}
2\inf(-h_{i}')u_{i} & > & \sum_{j=1,j\neq i}^{N}(\sup(\frac{\partial g_{i}}{\partial x_{j}}))u_{i}\\
 & & + \sup(\frac{\partial g_{j}}{\partial x_{i}})u_{j})
\end{array}\right.$$
for $i = 1,...,N$, then the Frechet derivative of (\ref{eq:1}) at $(u_{1},u_{2},...,u_{N})$ is invertible, where the $inf$ and $sup$ are defined on $R$ and $R^{n}$, respectively.
\end{thm} 
\pf The solutions operator for (\ref{eq:1}) is $A : (C^{2,\alpha})^{N} \rightarrow (C^{\alpha})^{N}$ such that for all $(v_{1},...,v_{N}) \in (C^{2,\alpha})^{N}$,
$$A((v_{1},...,v_{N})) = (w_{1},...,w_{N}),$$
where $w_{i} = -\Delta v_{i} - v_{i}(h_{i}(v_{i}) - g_{i}(v_{1},...,v_{i-1},v_{i+1},...,v_{N}))$ for $i = 1,...,N.$   
The Frechet derivative of $A$ at $(u_{1},...,u_{N})$ is $B = (a_{ij})$, where\\
$$\left.\begin{array}{l}
a_{ii} = -\Delta-(h_{i}(u_{i})-g_{i}(u_{1},...,u_{i-1},u_{i+1},...,u_{N}))-u_{i}h_{i}'(u_{i}),\\
a_{ij} = u_{i}\frac{\partial g_{i}(u_{1},...,u_{i-1},u_{i+1},...,u_{N})}{\partial u_{j}}
\end{array}\right.$$ 
for $i,j = 1,...,N, i\neq j.$\\
We need to show that $N(B) = \{0\}$ by Fredholm alternative. If
$$\left.\begin{array}{lll}
-\Delta \varphi_{i} & - & (h_{i}(u_{i}) - g_{i}(u_{1},...,u_{i-1},u_{i+1},...,u_{N}) + u_{i}h_{i}'(u_{i}))\varphi_{i}\\
& & + u_{i}\frac{\partial g_{i}(u_{1},...,u_{i-1},u_{i+1},...,u_{N})}{\partial u_{1}}\varphi_{1}
+ ... + u_{i}\frac{\partial g_{i}(u_{1},...,u_{i-1},u_{i+1},...,u_{N})}{\partial u_{i-1}}\varphi_{i-1}\\
& & + u_{i}\frac{\partial g_{i}(u_{1},...,u_{i-1},u_{i+1},...,u_{N})}{\partial u_{i+1}}\varphi_{i+1} + ... + u_{i}\frac{\partial g_{i}(u_{1},...,u_{i-1},u_{i+1},...,u_{N})}{\partial u_{N}}\varphi_{N} = 0,
\end{array}\right.$$
for $i = 1,...,N,$
then 
$$\left. \begin{array}{lll}
\int_{\Omega}[|\nabla \varphi_{i}|^{2} & - & (h_{i}(u_{i}) - g_{i}(u_{1},...,u_{i-1},u_{i+1},...,u_{N}) + u_{i}h_{i}'(u_{i}))\varphi_{i}^{2}\\
& + & (\frac{\partial g_{i}(u_{1},...,u_{i-1},u_{i+1},...,u_{N})}{\partial u_{1}}\varphi_{1}
 + ... + \frac{\partial g_{i}(u_{1},...,u_{i-1},u_{i+1},...,u_{N})}{\partial u_{i-1}}\varphi_{i-1}\\
& + & \frac{\partial g_{i}(u_{1},...,u_{i-1},u_{i+1},...,u_{N})}{\partial u_{i+1}}\varphi_{i+1} + ... + \frac{\partial g_{i}(u_{1},...,u_{i-1},u_{i+1},...,u_{N})}{\partial u_{N}}\varphi_{N})u_{i}\varphi_{i}] = 0,
\end{array}\right.$$
for $i = 1,...,N$.
Since $\lambda_{1}(g_{i}(u_{1},...,u_{i-1},u_{i+1},...,u_{N}) - h_{i}(u_{i})) = 0$ for $i = 1,...,N,$
$$\left. \begin{array}{l}
\int_{\Omega}[|\nabla \varphi_{i}|^{2} - (h_{i}(u_{i}) - g_{i}(u_{1},...,u_{i-1},u_{i+1},...,u_{N}))\varphi_{i}^{2} \geq 0 
\end{array} \right.$$     
for $i = 1,...,N$.
Hence,
$$\left. \begin{array}{l}
\int_{\Omega}-u_{i}h_{i}'(u_{i})\varphi_{i}^{2} + (\frac{\partial g_{i}(u_{1},...,u_{i-1},u_{i+1},...,u_{N})}{\partial u_{1}}\varphi_{1} + ... + \frac{\partial g_{i}(u_{1},...,u_{i-1},u_{i+1},...,u_{N})}{\partial u_{i-1}}\varphi_{i-1}\\
+ \frac{\partial g_{i}(u_{1},...,u_{i-1},u_{i+1},...,u_{N})}{\partial u_{i+1}}\varphi_{i+1} + ... + \frac{\partial g_{i}(u_{1},...,u_{i-1},u_{i+1},...,u_{N})}{\partial u_{N}}\varphi_{N})u_{i}\varphi_{i} \leq 0,
\end{array}\right.$$
for $i = 1,...,N.$
Hence, 
$$\int_{\Omega}\sum_{i=1}^{N}-u_{i}h_{i}'(u_{i})\varphi_{i}^{2} + \sum_{i=1}^{N}u_{i}\varphi_{i}\sum_{j=1,j\neq i}^{N}\frac{\partial g_{i}(u_{1},...,u_{i-1},u_{i+1},...,u_{N})}{\partial u_{j}}\varphi_{j} \leq 0.$$
It implies that
$$\int_{\Omega}\sum_{i=1}^{N}(-u_{i}h_{i}'(u_{i})\varphi_{i}^{2} + \sum_{j=1,j\neq i}^{N}\frac{\partial g_{i}(u_{1},...,u_{i-1},u_{i+1},...,u_{N})}{\partial u_{j}}u_{i}\varphi_{j}\varphi_{i}) \leq 0.$$
But,
$$\left.\begin{array}{lll}
& &\frac{\partial g_{i}(u_{1},...,u_{i-1},u_{i+1},...,u_{N})}{\partial u_{j}}u_{i}\varphi_{i}\varphi_{j}\\
& \leq & \frac{\partial g_{i}(u_{1},...,u_{i-1},u_{i+1},...,u_{N})}{\partial u_{j}}u_{i}(\frac{\varphi_{i}^{2}}{2} + \frac{\varphi_{j}^{2}}{2}).
\end{array}\right.$$
If
$$-u_{i}h_{i}'(u_{i}) > \sum_{j=1,j\neq i}^{N}(\frac{\frac{\partial g_{i}(u_{1},...,u_{i-1},u_{i+1},...,u_{N})}{\partial u_{j}}u_{i}}{2} + \frac{\frac{\partial g_{j}(u_{1},...,u_{j-1},u_{j+1},...,u_{N})}{\partial u_{i}}u_{j}}{2})$$
for $i = 1,...,N$, then the integrand in above inequality is positive definite, which implies $(\varphi_{1},...,\varphi_{N})$ is trivial. But, it holds if
$$\left.\begin{array}{lll}
2\inf(-h_{i}')u_{i} & > & \sum_{j=1,j\neq i}^{N}(\sup(\frac{\partial g_{i}}{\partial x_{j}})u_{i}\\
 & & + \sup(\frac{\partial g_{j}}{\partial x_{i}})u_{j})
\end{array}\right.$$
for $i = 1,...,N$, where the $inf$ and $sup$ are defined on $R$ and $R^{n}$, respectively.
\\[0.2in]
Combining the Theorems \ref{thm:1.1}, \ref{thm:3.1} and \ref{thm:3.3}, we have the following which is actually the main result in this section.
\begin{cor}\label{cor:3.4} Suppose\\
$(A)$ $h_{i}(0) > \lambda_{1} + g_{i}(k_{1},...,k_{i-1},k_{i+1},...,k_{N})$ for $i = 1,...,N$ and\\
$(B)$ 
$$\left.\begin{array}{lll}
-2\sup(h_{i}') & > & \sum_{j=1,j\neq i}^{N}[\sup(\frac{\partial g_{i}}{\partial x_{j}}) + K\sup(\frac{\partial g_{j}}{\partial x_{i}})],
\end{array}\right.$$
where $K = \sup_{\Omega,i\neq j}\frac{\theta_{h_{j}}}{\theta_{h_{i} - g_{i}(k_{1},...,k_{i-1},k_{i+1},...,k_{N})}}$.\\
Then there is a neighborhood $V$ of $(h_{1},...,h_{N})$ in $(C_{B}^{1,\alpha})^{N}$ such that if\\
$(\bar{h_{1}},...,\bar{h_{N}}) \in V$, then 
(\ref{eq:1}) with $(h_{1},...,h_{N}) = (\bar{h_{1}},...,\bar{h_{N}})$ has a unique coexistence state.  
\end{cor}
\pf From $\theta_{h_{i}} < k_{i}$ and the monotonicity of $g_{i}$ for $i = 1,...,N$, we have
$$\left\{ \begin{array}{lll}
h_{i}(0) & > & \lambda_{1} + g_{i}(k_{1},...,k_{i-1}k_{i+1},...,k_{N})\\ 
 & \geq & \lambda_{1}(g_{i}(\theta_{h_{1}},...,\theta_{h_{i-1}},\theta_{h_{i+1}},...,\theta_{h_{N}})). 
\end{array} \right.$$
The condition already guarantees that there is a unique coexistence state $(u_{1},...,u_{N})$ from Theorem \ref{thm:1.1}.\\
Furthermore, by the definition of $K$ and the estimate of the solution in the proof of Theorem \ref{thm:1.1}, we obtain  
$$2\inf(-h_{i}')u_{i} > \sum_{j=1,j\neq i}^{N}[\sup(\frac{\partial g_{i}}{\partial x_{j}})u_{i} + \sup(\frac{\partial g_{j}}{\partial x_{i}})u_{j}.$$ 
It implies that the Frechet derivative of (\ref{eq:1}) at $(u,v)$ is invertible from Theorem \ref{thm:3.3}. Therefore, the theorem follows from Theorem \ref{thm:3.1}.

\end{document}